\documentclass{article}

\usepackage{amssymb}
\usepackage{amsmath}
\usepackage{amsfonts}
\usepackage{float}
\usepackage{epsfig}
\usepackage{graphicx}

\newtheorem{remark}{Remark}

\def\vecz{\mathbf{z}}
\def\imi{\mathbf{i}}

\begin{document}

\title{Time-Limited and k-Limited Polling Systems: \\ A Matrix Analytic Solution}

\author{Ahmad Al Hanbali, Roland de Haan, Richard J. Boucherie,
\\ and Jan-Kees van Ommeren\\
University of Twente, Enschede, The Netherlands
}

\maketitle

\begin{abstract}
  In this paper, we will develop a tool to analyze polling systems with the
  autonomous-server, the time-limited, and the k-limited service discipline.
  It is known that these disciplines do not satisfy the well-known branching
  property in polling system, therefore, hardly any exact result exists in the literature
  for them. Our strategy is to apply an iterative scheme that
  is based on relating in closed-form the joint queue-length at the beginning
  and the end of a server visit to a queue. These kernel relations are derived
  using the theory of absorbing Markov chains. Finally, we will show that our tool
  works also in the case of a tandem queueing network with a single server that
  can serve one queue at a time.\\
  \newline


\end{abstract}

%\category{C.4}{Performance of Systems}{Performance Attributes}
%\category{G.3}{Probability and Statistics}{Queueing Theory}
%A category including the fourth, optional field follows...
%\keywords{Absorbing Markov chains; Matrix analytic solution; Polling system;
%   Autonomous-server discipline; Time-limited discipline; $k$-limited discipline;
%   Iterative scheme; Performance analysis;}

\section{Introduction}
\label{sec:intro}

Polling systems have been extensively studied in the last years due
to their vast area of applications in production and telecommunication
systems~\cite{sidi,takagi_2000}.  They have demonstrated to offer an
adequate modeling framework to analyze systems in which a set of
entities need certain service from a single resource.  These entities
are located at different positions in the system awaiting their turn
to receive service.

In queueing theory, a polling system is equivalent to a set of queues
with exogenous job arrivals all requiring an amount of service from a
single server.  The server serves each queue according to a specific
service discipline and after serving a queue he will move to a next
queue.  A key role in the analysis of such polling systems is played
by the so-called branching property~\cite{resing}.  This property
states that each job present at a queue at the arrival instant of the
server will be replaced in an independent and indentically distributed manner by a random number of
jobs during the course of the server's visit.  Service disciplines satisfying
the branching property yield a tractable analysis, while for
disciplines not satisfying this property hardly any exact results are
known.

The two most well-known disciplines that satisfy the bran- ching
property are the exhaustive and gated discipline.  Exhaustive means
that the server continues servicing a queue until it becomes empty.
At this instant the server moves to the next queue in his schedule.
Gated means that the server only serves the jobs present in the queue
at its arrival.

The drawback of the exhaustive and gated disciplines is that the
server is controlled by the job arrivals.  To reduce this control on
the server, other type of service disciplines were introduced such as
the time-limited and the $k$-limited discipline.
%However, we should emphasize that these policies do not satisfy the branching property and thus are intrinsically hard to analyze.
%These two policies verify the so-called branching property which makes their
%analysis tractable~\cite{resing}. This is because this property allows
%to simply and directly relate the joint number of jobs at the
%beginning and end of a visit period to a certain queue.
%The drawback of
%the exhaustive and gated disciplines is that the server is controlled
%by the job arrivals. For this reason, other type of polling policies
%were introduced to reduce the control on the server such as the time-limited and the $k$-limited discipline.
According to the time-limited discipline, the server continues
servicing a queue for a certain time period or until the queue becomes empty,
whichever occurs first.
%When the server is active at the end of a server
%visit, service will be either preempted or not.
%assumed. In the preemptive case the server leaves immediately after the
%timer expiration however the server leaves when it finishes servicing this
%job.
Under the k-limited discipline, the server continues servicing a queue
until $k$ jobs are served or the queue becomes empty, whichever occurs
first.  Another discipline, evaluated more recently in the literature
and closely related to the time-limited discipline, is the so-called
autonomous-server discipline~\cite{Asmta,Haan} which works as follows.
The server continues servicing a queue for a certain period of time
despite that, meanwhile, the queue may become empty.  This discipline
may also be seen as the non-exhaustive time-limited discipline.  We
should emphasize that these latter disciplines do not verify the
branching property and thus hardly any closed-form results are known
for the queue-length distribution under these disciplines.

To circumvent this difficulty, researchers resort to numerical methods
using for instance iterative solution techniques
%Discrete Fourier Transforms (DFT)
or by using a power series algorithm.  The power series algorithm
\cite{Blanc1,Blanc2} aims at solving the global balance equations.  To
this end, the state probabilities are written as a power series and
via a complex computation scheme the coefficients of these series, and
thus the queue-length probabilities, are obtained.  The iterative
techniques \cite{Leung1,Leung2} exploit the relations between the
joint queue-length distributions at specific instants, viz., the start
of a server visit and the end of a server visit.  The relation between
the queue length at the start and end of a visit to a queue is
established via recursively expressing the queue length at a job
departure instant in terms of the queue length at the previous
departure instant of a job.  The complementary relation, between the
queue length at the end of a visit to a queue and a start of visit to
a next queue, can easily be established via the switch-over time.
Starting with an initial distribution, the stationary queue-length
distribution is then obtained by means of iteration.  Although these
methods offer a way to numerically solve intrinsically hard systems,
their solution provides little fundamental insight and moreover the
computation time and memory requirements to obtain this solution are
exponential functions of the number of queues.

% such as the discrete
%Fourier transforms technique along with an iterative scheme for time-limited
%discipline~\cite{Leung2}.
%The main disadvantage of these numerical approaches is that time and
%memory requirements are exponential functions of the number of queues.
%
%This is due to the fact that to relate the number of jobs at the beginning and the
%end of a visit period to a queue the author in~\cite{Leung1,Leung2} conditions
%on all possible events of job departure during that visit period ???????.

In this paper, we develop a tool to analyze the autonomous server, the
time-limited, and the k-limited discipline. Our tool incorporates an
iterative solution method which enhances the method introduced
in~\cite{Leung1}.  More specifically, contrary to that approach, we will
establish a direct and more insightful relation between the joint number of
jobs at the beginning and end of a visit period to a queue
without conditioning on any intermediate events that occur during a
visit.  To this end,
%To relate the joint number of jobs at the beginning and end
we use the theory of absorbing Markov chains (AMC)~\cite{absorbing,neuts}.
We construct an AMC whose transient states represent the states of the polling system.
The event of the server leaving a queue is modeled as an absorbing event.
%Moreover, the states of this AMC represent the states of the polling system.
We will set the initial state of the AMC to the joint number of jobs
at the beginning of a service period of a queue.  Therefore, to find
the joint number of jobs at the end of a service period, it is
sufficient to keep track of the state from which the transition to
the absorption state occurs. The probability of the latter event is eventually
determined by first ordering the states in a careful way and
consequently exploiting the structures that arise in the generator matrix
of the AMC.  Following this approach, we relate in closed-form the joint
queue-length probability generating functions (p.g.f.)\ at the end of a
visit period to a queue to the joint queue-length p.g.f.\ at the
beginning of this visit period.  The major part of this paper is
devoted to deriving these kernel relations for the above-mentioned
three disciplines: autonomous-server, time-limited, and k-limited.  Once
these relations are obtained, the joint queue-length distribution at
server departure instants is readily obtained via a simple iterative
scheme.

The paper is organized as follows. In Section \ref{sec:model} we give a
careful description of the model and the assumptions. Section \ref{sec:auto}
analyses the autonomous-server discipline. In Section \ref{sec:time-limit} we
study the time-limited discipline. Section \ref{sec:k-limit} evaluates the
k-limited discipline. In Section~\ref{sec:ite} we describe the iterative
scheme that is important to compute the joint queue-length distribution.
Section \ref{sec:tand} analyses briefly the tandem model case with the
autonomous-server and the time-limited service discipline. Finally, in
Section \ref{sec:conc}, we conclude the paper and give some research directions.
\section{Model}
\label{sec:model}

We consider a single-server polling model consisting of $M$
first-in-first-out (FIFO) systems with unlimited queue, $Q_i$,
$i=1,\ldots,M$. Jobs arrive to $Q_i$ according to a Poisson process
with arrival rate $\lambda_i$.  We let $N_i(t)$ denote the number of
jobs in $Q_i$, $i=1,\ldots,M$, at time $t \geq 0$ and it is assumed
that $N_i(0)=0$, $i=1,\ldots,M$.  The service requirement $B_i$ at
$Q_i$ has an exponential distribution $B_i(\cdot)$ and mean $b_i$. We
assume that the service requirements are independent and identically
distributed (iid) random variables (rvs).  The server visits the queues in a cyclic
fashion.  After a visit to $Q_i$, the server incurs a switch-over time
$C^{i}$ from $Q_i$ to $Q_{i+1}$.  We assume that $C^{i}$ is
independent of the service requirement and follows a general
distribution $C^{i}(\cdot)$ with mean $c^{i}$, where at least one
$c^i>0$. The service discipline at each queue is either
autonomous-server, time-limited, or $k$-limited. It is assumed that
the queues of the polling system are stable.

%Under the autonomous-server discipline, the server remains at location $Q_i$ an exponentially distributed time
%with rate $\alpha_i$ before it migrates to the next queue in the cycle.
%It is stressed that even when $Q_i$ becomes empty, the server will remain at this queue.
%Under the time-limited discipline, the server departs from $Q_i$ when it becomes empty or when an exponentially distributed time
%with rate $\alpha_i$ has elapsed, whichever occurs first.
%Finally, under the k-limited discipline, the server departs from $Q_i$ when it has become empty or when $k$ jobs are served,
%whichever occurs first.

In case the server is active at the end of a server visit, which may
happen under the autonomous-server and time-limited disciplines, then
the service will be preempted.  At the beginning of the next visit of
the server, the service time will be re-sampled according to
$B_i(\cdot)$. This discipline is commonly referred to as {\it preemptive-repeat-random}.

A word on notation. Given a random variable $X$, $X(t)$ will denote
its distribution function.  We use {\bf I} to denote an identity matrix of
appropriate size and use $\otimes$ as tensor product operator defined as follows.
Let {\bf A} and {\bf B} be two matrices and $a(i,j)$ and
$b(i,j)$ denote the (i,j)-entries of {\bf A} and {\bf B} respectively
then ${\bf A}\otimes {\bf B}$ is a block matrix where the (i,j)-block is equal to
b(i,j){\bf A}. We use $e$ to denote a row vector of elements equal to one
and $e_i$ to denote a row vector with the $i$-th element equal to one and the
other elements equal to zero. Finally, $v^T$ will denote the transpose of vector $v$.

% $\otimes$ is the Kronecker product operator.
%%% $\forall$ means for all elements.

\section{Autonomous-server discipline}
\label{sec:auto}

In this section, we will relate the joint queue-length probabilities
at the beginning and end of a server visit to a queue for the
autonomous-server discipline.  Under the autonomous-server discipline,
the server remains at location $Q_i$ an exponentially distributed time
with rate $\alpha_i$ before it migrates to the next queue in the
cycle.  It is stressed that even when $Q_i$ becomes empty, the server
will remain at this queue.

Without loss of generality let us consider a server visit to $Q_1$. We assume that
the p.g.f.~of the steady-state queue-length at service's beginning instant at $Q_1$,
denoted by $\beta^{A}_1({\bf z})$, is known, where ${\bf z}:=(z_1,\ldots,z_M)$ and $|z_i|\leq 1$ for $i=1,\ldots,M$. The aim
is to derive the p.g.f.~of the steady-state queue-length at service visit's end at $Q_1$,
denoted by $\gamma^A_1(\mathbf{z})$. In order to compute $\gamma^A_1(\mathbf{z})$, we first assume
that $Q_1$ has a limited length of $L-1$ jobs including the job in service. This queue is denoted by $Q_1^L$.
Later, we will let $L$ tend to infinity to get the desired results.

The probability that there are $(i_1,...,i_M)$ jobs in
$(Q_1,\ldots,$ $Q_M)$ at the beginning of a server visit to $Q_1$ is
denoted by $\mathbb{P}_L\big({\bf N}_1^b=(i_1,\ldots,i_M)\big)$. Similarly, the
probability that there are $(j_1,...,j_M)$ jobs in $(Q_1,\ldots,Q_M)$
at the end of a server visit to $Q_1$ is denoted by
$\mathbb{P}_L\big({\bf N}_1^e=(j_1,\ldots,j_M)~|~{\bf N}_1^b=(i_1,\ldots,i_M)\big)$. Under the assumption that the
unlimited $Q_1$ is stable, $\lim_{L\to\infty}$ $\mathbb{P}_L\big({\bf N}_1^b=(i_1,\ldots,i_M)\big)=$
$\mathbb{P}\big({\bf N}_1^b=(i_1,\ldots,i_M)\big)$ and $\beta^{A}_1(\mathbf{z})=\mathbb{E}[{\bf z}^{{\bf N}_1^b}]$ .

Let ${\bf N}(t):=(N_1(t),\ldots,N_M(t))$ denote the $M$-dimensional,
continuous-time Markov chain with discrete state-space
$\xi_A=\{0,1,$ $\ldots,L-1\}\times\{0,1,\ldots\}^{M-1}\cup\{a\}$, where
$N_j(t)$ represents the number of jobs in $Q_j$ at time $t$. State $\{a\}$ is
absorbing. We refer to this absorbing Markov chain by ${\bf AMC}_A$. The absorption of ${\bf AMC}_A$
occurs when the server leaves $Q_1$ which happens with rate $\alpha_1$. Moreover,
the initial state of ${\bf AMC}_A$ at $t=0$ is set to the system state at server's
arrival to $Q_1$, i.e., $N_1^b=(i_1,\ldots,i_M)$. Therefore, the probability
that the absorption of ${\bf AMC}_A$ occurs from one of the states
$\{(j_1,\ldots,j_M)\}$ equals $\mathbb{P}_L\big({\bf N}_1^e=(j_1,\ldots,j_M)~|~N_1^b=(i_1,\ldots,i_M)\big)$. Let
${\bf n}=(n_1,\dots,n_M)\in \xi_A-\{a\}$ and $e_l$ the $M$-dimensional row vector whose
entries equal zero except the $l$-th entry that equals one. The non-zero
transition rates of ${\bf AMC}_A$ can be written as
\begin{eqnarray*}
\label{eq:transprob}
\begin{array}{lll}
q({\bf n},{\bf n}+e_1)          &=& \lambda_1,\quad~ 0 \leq n_1\leq L-2,\\
q({\bf n},{\bf n}+e_l)          &=& \lambda_l,\quad~ 2\leq l\leq M,\\
q({\bf n},{\bf n}-e_1)          &=& 1/b_1,~ 1 \leq n_1 \leq L-1,\\
q\big({\bf n},\{a\}\big)        &=& \alpha_1 .
\end{array}
\end{eqnarray*}

%\begin{eqnarray*}
%\begin{array}{lll}
%\end{array}
%\end{eqnarray*}

%\begin{figure}[hbtp]
%\vspace{0pt}\vspace{-0.4cm}
%  \centerline{\hbox{
%    \includegraphics[width=2.5in]{mc2}
%   }
%  }
%\vspace{0cm}
%\vspace{-0.2cm}
%\caption{Transition state diagram of ${\bf AMC}$.}
%  \label{fig:mc2}
%\vspace{-0.4cm}
%\end{figure}

We derive now $\mathbb{P}_L\big({\bf N}_1^e=(j_1,\ldots,j_M)|~{\bf N}_1^b=(i_1,\ldots,i_M)\big)$.
During a server visit to $Q_1$, the number of jobs at $Q_l$, $l=$ $2,\ldots,M$, may only increase.
Therefore $\mathbb{P}_L\big({\bf N}_1^e=(j_1,\ldots,j_M$ $)~|~{\bf N}_1^b=(i_1,\ldots,i_M)\big)$ is
strictly positive for  $j_l\geq i_l$, $l=2,\ldots,M$, and zero otherwise.
For sake of clarity, we will show first in detail the structure of ${\bf AMC}_A$ in the case of
$3$ queues, i.e.~for $M=3$, before considering the general case.
\\ \newline
{\bf Case M=3}. Let us consider the transient states of ${\bf AMC}_A$,
i.e., $(n_1,n_2,n_3)\in \xi_A-\{a\}$, where
$n_1\in\{0,1,\dots,L-1\}$ and $n_2,n_3\in\{0,1,\dots\}$. We recall
that we consider a server visit to $Q_1$. The number of jobs at
$Q_2$ and $Q_3$ may only increase during a server visit to $Q_1$, while
the number of jobs at $Q_1$ may increase or decrease. To take advantage
of this property, we will order the transient states of the ${\bf AMC}_A$
as follows: $(0,0,0),(1,0,0),(2,0,0),\ldots,$
$(0,1,0),(1,1,0),(2,1,0),$ $\ldots,(0,0,1),(1,0,1),(2,0,1),\ldots$, i.e.,
lexicographically orde- red first according to $n_3$, then $n_2$, and finally
according to $n_1$. This ordering induces that the generator matrix of the transition rates between
the transient states of ${\bf AMC}_A$ for $M=3$, denoted by ${\bf Q}_3$,
satisfies the following structure. That is, ${\bf Q}_3$ is an infinite
upper-bidiagonal block matrix with diagonal blocks equal to ${\bf A}_3$ and
upper-diagonal blocks equal $\lambda_3{\bf I}$, i.e.,
\begin{equation}
{\bf Q}_3=\left(\begin{array}{lllll}
{\bf A}_3 & \lambda_3{\bf I} & {\bf 0}          & \cdots & \cdots\\
{\bf 0}   & {\bf A}_3        & \lambda_3{\bf I} & ~~{\bf 0} & \cdots \\
\vdots\ddots  &~~~~\ddots     &            & \ddots & \ddots\\
\end{array}\right).
\end{equation}
We note that ${\bf A}_3$ denotes the generator matrix of the
transitions which do not induce any modification in the number
of jobs at $Q_3$. Moreover, $\lambda_3{\bf I}$ denotes the
transition rate matrix between the transient states $(n_1,n_2,n_3)$
and $(n_1,n_2,n_3+1)$, i.e., the transitions that represent an arrival to $Q_3$.
The block matrix ${\bf A}_3$ is also an infinite
upper-bidiagonal block  matrix with diagonal blocks equal to
${\bf A}_2$, and upper-diagonal blocks equal $\lambda_2{\bf I}$, i.e.,
\begin{equation}
{\bf A}_3=\left(\begin{array}{lllll}
{\bf A}_2 & \lambda_2{\bf I} & {\bf 0}          & \cdots & \cdots\\
{\bf 0}   & {\bf A}_2        & \lambda_2{\bf I} & ~~{\bf 0} & \cdots \\
\vdots\ddots  &~~~~\ddots     &            & \ddots & \ddots\\
\end{array}\right),
\end{equation}
where $\lambda_2{\bf I}$ denotes the transition rate matrix between
the transient states $(n_1,n_2,n_3)$ and $(n_1,n_2+1,n_3)$ and
${\bf A}_2$ is the generator matrix of the transition between
the transient states $(n_1,n_2,n_3)$ and $(n_1\pm 1,n_2,n_3)$.
Observe that ${\bf A}_2$ equals the sum of the generator
matrix of an M/M/1/L-1 queue with arrival rate $\lambda_1$ and departure
rate $1/b_1$ and of the matrix $-(\lambda_2+\lambda_3+\alpha_1)${\bf I}.
Now, we compute $\mathbb{P}_L\big({\bf N}_1^e=(j_1,j_2,j_3)~|~{\bf N}_1^b=(i_1,i_2,i_3)\big)$
as function of the inverse of  ${\bf Q}_3$, ${\bf A}_3$ and ${\bf A}_2$. First note that
since ${\bf Q}_3$, ${\bf A}_3$ and ${\bf A}_2$ are all sub-generators with sum of their row
elements strictly negative, these matrices are invertible. From the
theory of absorbing Markov chains, given that {\bf AMC}$_A$ starts in state $(i_1,i_2,i_3)$,
the probability that the transition to the absorption state $\{a\}$ occurs from state $(j_1,j_2,j_3)$
reads (see, e.g., \cite{gaver})
\begin{equation}
\mathbb{P}_L\big({\bf N}_1^e=(j_1,j_2,j_3)~|~{\bf N}_1^b\big)=-\alpha_1c_3({\bf Q}_3)^{-1}d_3,
\label{eq:cprobn1e}
\end{equation}
where $c_3$ is the probability distribution vector of {\bf AMC}$_A$'s initial state that can be given
by
\[
c_3:=e_{i_1} \otimes e_{i_2} \otimes e_{i_3},
\]
and $\alpha_1d_3$ is the transition rate vector to $\{a\}$ given that $(j_1,j_2,j_3)$ is the last state
visited before absorption where $d_3$ can be given by
\begin{equation*}
d_3:=(e_{j_1} \otimes e_{j_2} \otimes e_{j_3})^T.
\end{equation*}
${\bf Q}_3$ is an upper-bidiagonal block matrix. Hence, it is easy to show
that $({\bf Q}_3)^{-1}$ is an upper-triangular block matrix with (i,j)-block equal to
$(-({\bf A}_3)^{-1}\lambda_3{\bf I})^{j-i}({\bf A}_3)^{-1}$, thus we find that
\begin{eqnarray}
c_3({\bf Q}_3)^{-1}d_3&=&c_2(-\lambda_3({\bf A_3})^{-1})^{j_3-i_3}({\bf A}_3)^{-1}d_2,\label{eq:c3q3d3}
\end{eqnarray}
where $c_2=e_{i_1} \otimes e_{i_2}$ and $d_2= (e_{j_1} \otimes e_{j_2})^T$. Plugging (\ref{eq:c3q3d3})
into (\ref{eq:cprobn1e}) gives that
\begin{eqnarray}
\mathbb{P}_L\big({\bf N}_1^e=(j_1,j_2,j_3)~|~{\bf N}_1^b\big)=-\alpha_1c_2(-\lambda_3({\bf A_3})^{-1})^{j_3-i_3}({\bf A}_3)^{-1}d_2.
\end{eqnarray}
\newline
{\bf General case}. By analogy with the case of $M=3$, we order the transient states of
${\bf AMC}_A$ first according to $n_M$, then $n_{M-1}$, $\ldots$, and finally
according to $n_1$. During a server visit to $Q_1$, the number of jobs at
$Q_j$, $j=2,\ldots,M$, may only increase. Therefore, similarly to the case of $M=3$,
the {\bf AMC}$_A$ the generator matrix of the transition rates between
the transient states of ${\bf AMC}_A$ for the general case, denoted by
${\bf Q}_M$, is an upper-bidiagonal block matrix with diagonal blocks equal to
${\bf A}_M$, and upper-diagonal blocks equal to $\lambda_M{\bf I}$. Moreover, ${\bf A}_M$ in turn
is an upper-bidiagonal block matrix with diagonal blocks equal to ${\bf A}_{M-1}$,
and upper-diagonal blocks equal to $\lambda_{M-1}{\bf I}$.  We emphasize that
${\bf A}_j$, $j=M,\ldots,3$, all verify the previous property.  Finally, the matrix ${\bf A}_2$
equals the sum of the generator matrix of an M/M/1/L-1 queue with
arrival rate $\lambda_1$ and departure rate $1/b_1$ and of the matrix
$-(\lambda_2+\ldots+\lambda_M+\alpha_1)${\bf I}.

By analogy with the $M=3$ case, we find that the probability of ${\bf N}_i^e=(j_1,\ldots,j_M)$,
given that ${\bf N}_1^b=(i_1,\ldots,i_M)$, reads
\begin{eqnarray}
\mathbb{P}_L\big({\bf N}_1^e=(j_1,\ldots,j_M)~|~{\bf N}_1^b=(i_1,\ldots,i_M)\big)= \nonumber ~~~~~~~~~~~~~~~~~~~~~~~~~~~~~~~~~~~~~~~~&&\\
-\alpha_1 c_{M-1}\big(-\lambda_M({\bf A_M})^{-1}\big)^{j_M-i_M}({\bf A_M})^{-1}d_{M-1}.&&
\label{eq:conprob}
\end{eqnarray}
\begin{eqnarray}
c_{M-1}&:=& e_{i_1}\otimes\ldots\otimes e_{i_{M-1}},\nonumber\\
d_{M-1}&:=&(e_{j_1}\otimes\ldots\otimes e_{j_{M-1}})^T.\nonumber
\label{eq:xx1}
\end{eqnarray}
We derive now the conditional p.g.f.~of ${\bf N}_1^e$. Note that
$\big(-\lambda_M({\bf A}_M)^{-1}\big)$ is a sub-stochastic matrix
with the sum of its row elements strictly smaller than one,
which gives that $\lim_{n\to\infty}(-\lambda_M({\bf A}_M)^{-1})^n={\bf 0}$. Combining the
latter result with (\ref{eq:conprob}) we find that
\begin{eqnarray}
\mathbb{E}_L\Big[{\bf z}^{{\bf N}^e_1}~|~{\bf N}_1^b=(i_1,\ldots,i_M) \Big]=
-\alpha_1 z_M^{i_M}c_{M-1}\big({\bf A}_M+z_M\lambda_M{\bf I}\big)^{-1}d_{M-1}({\bf z}),&&
\label{eq:elm}
\end{eqnarray}
where
\begin{eqnarray}
d_{M-1}({\bf z}):=
\sum_{j_1=0}^{L-1}\sum_{j_2\geq i_2}\ldots\sum_{j_{M-1}\geq i_{M-1}}(z_1^{j_1} e_{j_1}\otimes\ldots\otimes z_{M-1}^{j_{M-1}}e_{j_{M-1}})^T,&
\label{eq:xxx2}
\end{eqnarray}
and $|z_i|\leq 1$, $i=1,\ldots,M$. It remains to find $({\bf A}_M+z_M\lambda_M{\bf I})^{-1}$. Since
${\bf A}_M$ is an upper-bidiagonal block matrix, the (i,j)-block of
$({\bf A}_M+z_M\lambda_M{\bf I})^{-1}$ is given by $(-\lambda_{M-1})^{j-i}\times$ $({\bf A}_{M-1}+z_M\lambda_M{\bf I})^{-j+i-1}$.
Plugging the latter result into (\ref{eq:elm}) gives that
\begin{eqnarray}
E_L\Big[{\bf z}^{{\bf N}^e_1}~|~{\bf N}_1^b=(i_1,\ldots,i_M) \Big]=-\alpha_1 z_M^{i_M} z_{M-1}^{i_{M-1}}c_{M-2}\times ~~~~~~~~~~~~~~~~~~~~~~~~~&&\nonumber\\
\big({\bf A}_{M-1}+(z_M\lambda_M+z_{M-1}\lambda_{M-1}){\bf I}\big)^{-1}d_{M-2}({\bf z}),&&
\label{eq:xxx3}
\end{eqnarray}
 where
\begin{eqnarray}
c_{M-2}&:=& e_{i_1}\otimes\ldots\otimes e_{i_{M-2}},\nonumber\\
d_{M-2}({\bf z})&:=&
\sum_{j_1=0}^{L-1}\sum_{j_2\geq i_2}\ldots\sum_{j_{M-2}\geq i_{M-2}}(z_1^{j_1} e_{j_1}\otimes\ldots\otimes z_{M-2}^{j_{M-2}}e_{j_{M-2}})^T.\nonumber
\label{eq:xxx4}
\end{eqnarray}
By an induction argument along with the properties that ${\bf A}_j$, $j=3,\ldots,M-1$, is an upper-bidiagonal block matrix,
it can be shown that
\begin{eqnarray}
\mathbb{E}_L\Big[{\bf z}^{{\bf N}^e_1}~|~{\bf N}_1^b=(i_1,\ldots,i_M) \Big]=-\alpha_1 z_{2}^{i_{2}} \ldots z_M^{i_M} e_{i_1} \times ~~~~~~~~~~~~~~~~~~~~~~~~~~~~&& \nonumber \\
 \Big({\bf A}_{2}+(z_{2}\lambda_{2}+\ldots+z_M\lambda_M){\bf I}\Big)^{-1}d_{1}(z_1),&&
\label{eq:ztran0}
\end{eqnarray}
where
\begin{eqnarray}
d_{1}(z_1)&:=&\sum_{j_1=0}^{L-1}z_1^{j_1}(e_{j_1})^T=(1,z_1,\ldots,z_1^{L-1})^T.\nonumber
\end{eqnarray}
Removing the condition on ${\bf N}_1^b$, it is readily seen that
\begin{equation}
\mathbb{E}_L\big[{\bf z}^{{\bf N}^e_1}\big]=-\alpha_1 f \Big({\bf A}_{2}+(z_{2}\lambda_{2}+\ldots+z_M\lambda_M){\bf I}\Big)^{-1}d_{1}(z_1),
\label{eq:ztran1}
\end{equation}
where $f$ is the $L$-dimensional row vector with $i$-th element equal to
$\mathbb{E}\big[{\bf 1}_{\{N_1^b=i\}} \cdot z_2^{N_2^b}\ldots z_M^{N_M^b}\big]$, for $i=0,\ldots,L-1$.
It remains to find the inverse of ${\bf A}_{2}+(z_{2}\lambda_{2}+\ldots+z_M\lambda_M){\bf I}$
and to let $L\to\infty$.

Let $u^T=(1,0,\ldots,0)$ and let $v^T=(0,\ldots,0,1)$. We recall that
${\bf A}_2$ equals the sum of the generator matrix of an M/M/1/L-1 queue
with arrival rate $\lambda_1$ and departure rate $1/b_1$ and of the matrix
$-(\lambda_2+\ldots+\lambda_M+\alpha_1)${\bf I}.  Let
${\bf Q_{A}(z)}:={\bf A}_{2}+(z_{2}\lambda_{2}+\ldots+z_M\lambda_M){\bf I}$.
Now, observe that
${\bf Q_{A}(z)}={\bf T_{A}(z)}+1/b_1 uu^T+\lambda_1 vv^T$, where ${\bf T_{A}(z)}$
is a L-by-L tridiagonal Toeplitz matrix with diagonal entries equal
$\big(-\lambda_1-1/b_1-\alpha_1-\sum_{m=2}^{M}\lambda_m(1-z_m)\big)$,
upper-diagonal entries equal $\lambda_1$, and lower-diagonal
entries $1/b_1$.  Let $t_{ij}^*$ denote the $(i,j)$-entry of
${\bf T_A^{-1}(z)}$.  By applying the Sherman-Morrison formula \cite[p.
76]{Sherman-Morrison} we find that the $(i,j)$-entry of ${\bf
  Q_A^{-1}(z)}$ gives for $i,j=1,\ldots,L$,
\begin{equation}
q_{ij}^*=m_{ij}-\lambda_1\frac{m_{iL}m_{Lj}}{1+\lambda_1 m_{LL}},~~\mbox{where}~
 m_{ij}=t_{ij}^*-\frac{t_{i1}^*t_{1j}^*}{b_1+ t_{11}^*}.
\label{eq:invQz1}
\end{equation}
%with
%\begin{equation}
%m_{ij}=t_{ij}^*-\beta_1\frac{t_{i1}^*t_{1j}^*}{1+\beta_1 t_{11}^*},
%\label{eq:invQz2}
%\end{equation}
%for $i,j=1,\ldots,M$.

%Now we compute $t_{ij}$. Let $c_{-1}=z\beta_1$,
%$c_{0}=-(\lambda+\beta_1+\alpha_1)$, and $c_{1}=\lambda$. Let $r_{1,2}$
%denote the two roots of $c_{-1}+rc_0+r^2c_1$. Since $0\leq z\leq 1$,
%$\Delta^2:=c_0^2-4c_{-1}c_1>0$, thus these roots are distinct,
%real, and equal to
%\begin{equation}
%r_{1,2}=\frac{-c_0\mp\sqrt{c_0^2-4c_{-1}c_1}}{2c_1}.\label{eq:r1r2}
%\end{equation}
%Note that $r_{1,2}\geq 0$, $r_1<r_2$, $r_1<1$, and $r_2>1$. The
%(i,j)-entry of ${\bf T}(z)$ can be written as function of $r_{1,2}$
%as follows (Hint: See~\cite[Sec. 3.1]{Dow2003} )
The inverse of a tridiagonal Toeplitz matrix is known in closed-form
(see \cite[Sec.~3.1]{Dow2003})
\begin{equation}
t_{ij}^*=\left\{\begin{array}{l l}
-\frac{(r_{11}^i-r_{21}^i)(r_{11}^{L+1-j}-r_{21}^{L+1-j})}{\lambda_1(r_{11}-r_{21})(r_{11}^{L+1}-r_{21}^{L+1})}            &,~i\leq j\leq L\\
\frac{(r_{11}^{-j}-r_{21}^{-j})(r_{11}^{L+1}r_{21}^i-r_{21}^{L+1}r_{11}^i)}{\lambda_1(r_{11}-r_{21})(r_{11}^{L+1}-r_{21}^{L+1})} &,~j\leq i\leq L\\
\end{array} \right.
\label{eq:invT}
\end{equation}
where $r_{11}$ and $r_{21}$ are the distinct roots of
\begin{equation}
P_1(r):=\lambda_1 r^2-s_1r+1/b_1,
\end{equation}
where $s_1:=\lambda_1+1/b_1+\alpha_1+\sum_{m=2}^{M}\lambda_m(1-z_m)$.
%that read
%\begin{equation}
%r_{1,2}=\frac{(\lambda_1+1/b_1+\alpha_1+)\mp\sqrt{(\lambda_1+1/b_1+\alpha_1)^2-4\lambda_1/b_1}}{2\lambda_1}.\label{eq:r1r2}
%\end{equation}
We take $|r_{11}|<|r_{21}|$. Note that $|\lambda_1r^2+1/b_1|<|-s_1r|$
for every $|r|=1$, thus Rouch\'e's theorem gives that $P_1(r)$ has
exactly one root inside the disk of radius one for all $|z_i|\leq 1$
(see, e.g.,~\cite{rouche}). For this reason, we have that
$|r_{11}|<1<|r_{21}|$.

Inserting the values of $t_{ij}^*$ into~(\ref{eq:ztran1}) yields that
\begin{eqnarray}
\mathbb{E}_L\big[{\bf z}^{{\bf N}^e_1} \big]&=&-\alpha_1\sum_{i=0}^{L-1} f(i)\sum_{j=1}^{L}z_1^{j-1}\bigg[t_{ij}^*-\frac{1/b_1 t_{i1}^*t_{1j}^*}{1+1/b_1 t_{11}^*}
\nonumber\\
&&-\frac{\lambda_1 m_{iL}}{1+\lambda_1 m_{LL}}\bigg(t^*_{Lj}-\frac{1/b_1t_{L1}^*t_{1j}^*}{1+1/b_1 t_{11}^*}\bigg)\bigg].
\label{eq:ztran2}
\end{eqnarray}
Thus, it remains to let $L\to\infty$ in (\ref{eq:ztran2}) in order to find $\mathbb{E}\big[{\bf z}^{{\bf N}^e_1}\big]$.
It is readily seen that
\begin{eqnarray}
\lim_{L\to\infty}t^*_{LL-j}&=&-\frac{1}{\lambda_1 r_{21}}r_{11}^{j}, \nonumber \\
\lim_{L\to\infty}m_{L-iL}&=&\lim_{L\to\infty}t^*_{L-iL}=-\frac{1}{\lambda_1}r_{21}^{-(i+1)},\nonumber \\
\lim_{L\to\infty}t^*_{1j}&=&-\frac{1}{\lambda_1}r_{21}^{-j}, \nonumber \\
\lim_{L\to\infty}t^*_{i1}&=&\frac{-1}{\lambda_1 r_{11} r_{21}}r_{11}^i. \nonumber
\end{eqnarray}
Some technical calculus shows that the following limit is equal to zero
\begin{eqnarray}
\lim_{L\to\infty}\alpha_1\sum_{i=0}^{L-1} f(i)\sum_{j=1}^{L}z_1^{j-1}\frac{\lambda_1 m_{iL}}{1+\lambda_1 m_{LL}}\bigg(t_{Lj}-\frac{1/b_1 t_{L1}t_{1j}}{1+1/b_1 t_{11}}\bigg).\nonumber
\end{eqnarray}
Finally, plugging the previous limits in (\ref{eq:ztran2}) it can be shown that
\begin{equation}
%\mathbb{E}\big[{\bf z}^{{\bf N}^e_1}\big]&=&\frac{\alpha_1(1-z_1)}{P_1(z_1)}\bigg(\frac{r_{11}\mathbb{E}\big[({\bf z^*})^{{\bf N}^b_1}\big]}{1-r_{11}}-\frac{z_1\mathbb{E}[{\bf z}^{{\bf N}^b_1}]}{1-z_1}\bigg),\nonumber \\
\mathbb{E}\big[{\bf z}^{{\bf N}^e_1}\big]=\gamma^A_1({\bf z})=\frac{\alpha_1(1-z_1)}{P_1(z_1)}\bigg(\frac{r_{11}\beta^{A}_1({\bf z}^*_1)}{1-r_{11}}-\frac{z_1\beta^{A}_1({\bf z})}{1-z_1}\bigg),
\label{eq:ztransN}
\end{equation}
where ${\bf z}^*_1:=(r_{11},z_2,\ldots,z_M)$. Eq.~(\ref{eq:ztransN})
relates in closed-form $\gamma^A_1({\bf z})$, p.g.f.~of the joint
queue-length at the beginning of a server visit to $Q_1$, to
$\beta^{A}_1({\bf z})$, p.g.f.~of the joint queue-length at the end of a
server visit to $Q_1$. From (\ref{eq:ztransN}), we deduce that for a
server visit to $Q_i$, $i=1,\ldots,M$,
\begin{equation}
\gamma^A_i({\bf z})=\frac{\alpha_i(1-z_i)}{P_i(z_i)}\Big(\frac{r_{1i}\beta^{A}_i({\bf z}^*_i)}{1-r_{1i}}-\frac{z_i\beta^{A}_i({\bf z})}{1-z_i}\Big),
\end{equation}
where
\begin{eqnarray}
P_i(z_i) &:=& \lambda_i z_i^2-s_iz_i+1/b_i, \label{eq:pizi}\\
s_i      &:=& \lambda_i+1/b_i+\alpha_i+\sum_{m=1,m\neq i}^{M}\lambda_m(1-z_m), \label{eq:si}\\
r_{1i}   &:=& \frac{s_i-\sqrt{(s_i)^2-4\lambda_i/b_i}}{2\lambda_i}, \label{eq:r1i}
\end{eqnarray}
${\bf z}^*_i:=(z_1,\ldots,z_{i-1},r_{1i},z_{i+1},\ldots,z_M)$, and $|r_{1i}|<1$.

Finally, introducing the switch-over times from $Q_{i-1}$ to $Q_{i}$, thus by using that
$\beta^{A}_i({\bf z})=\gamma^A_{i-1}({\bf z})C^{i-1}({\bf z})$, where $C^{i-1}({\bf z})$ is the
p.g.f.~of the number of Poisson arrivals during $C^{i-1}$, we obtain
\begin{eqnarray}
\gamma^A_i({\bf z})&=&\frac{\alpha_i(1-z_i)r_{1i}}{P_i(z_i)(1-r_{1i})}\gamma^A_{i-1}({\bf z}^*_i)C^{i-1}({\bf z}^*_i)\nonumber\\
&&-\frac{\alpha_iz_i}{P_i(z_i)}\gamma^A_{i-1}({\bf z})C^{i-1}({\bf z}).
\end{eqnarray}

\section{Time-Limited discipline}
\label{sec:time-limit}
In this section, we will relate the joint queue-length probabilities at the beginning and end of a server visit to a queue
for the time-limited discipline.
Under this discipline, the server departs from $Q_i$ when it becomes empty or when a timer of exponential distribution
duration with rate $\alpha_i$ has expired, whichever occurs first.
Moreover, if the server arrives to an empty queue, he leaves the queue immediately
and jumps to the next queue in the schedule.
For this reason, we should differentiate here between the two events where the
server join an empty and non-empty queue. %As in Section~\ref{sec:auto}, our objective
%is to relate the joint queue-length p.g.f. $\gamma_1^{T}(\bold{z})$,~at the beginning
%and, $\beta_1^{T}(\bold{z})$, at end of a server visit to $Q_1$.

We will follow the same approach as in Section~\ref{sec:auto}. Thus,
we first assume that $Q_1$ has a limited queue of $L-1$ jobs, second there
are ${\bf N}_1^b:=(i_1,...,i_M)$ jobs in $(Q_1,\ldots,$ $Q_M)$, with $i_1\geq 1$, at the
beginning time of a server visit to $Q_1$ and third there are ${\bf N}_1^e:=(j_1,...,j_M)$ jobs
in $(Q_1,\ldots,Q_M)$ at the end time of a server visit to $Q_1$. Note that if $Q_1$
is empty at the beginning of a server visit, i.e.,~$i_1=0$,
$\mathbb{P}\big({\bf N}_1^e={\bf N}_1^b\big)=1$. We will exclude the latter
obvious case from the analysis in the following, however, we will include it
when we will uncondition on ${\bf N}_1^b$.

Let ${\bf N}(t):=(N_1(t),\ldots,N_M(t))$ denote the $M$-dimensional,
continuous-time Markov chain with discrete state-space
$\xi_{T}$ $=\{1,\ldots,L-1\}\times\{0,1,$ $\ldots\}^{M-1}\cup\{a\}$, where
$N_j(t)$ represents the number of jobs in $Q_j$ at time $t$ and
at which $Q_1$ is being served. State $\{a\}$ is absorbing. We refer to this
absorbing Markov chain by ${\bf AMC}_{T}$. The absorption of ${\bf AMC}_{T}$
occurs when the server leaves $Q_1$ which happens with rate $\alpha_1$ from all
transient states. The transient states of the form $(1,n_2,\ldots,n_M)$ have an additional transition
rate to $\{a\}$ that is equal to $1/b_1$, which represents the departure of the last job
at $Q_1$.

We set ${\bf N}(0)={\bf N}_1^b$. Therefore, the probability that the
absorption of ${\bf AMC}_{T}$ occurs from one of the states
$\{(j_1,\ldots,$ $j_M)\}$ equals $\mathbb{P}_L\big({\bf
  N}_1^e=(j_1,\ldots,j_M)\big)$, if the absorption is due to the timer
expiration with rate $\alpha_1$. However, if the absorption is due to
$Q_1$ becoming empty, $\mathbb{P}_L\big({\bf N}_1^e=(0,j_2\ldots,j_M)\big)$ equals
the probability that the absorption with rate $1/b_1$ occurs from one
of the states $\{(1,j_2,\ldots,j_M)\}$. The non-zero transition rates of
${\bf AMC}_{T}$ can be written for all {\bf n} $\in \xi_{T}-{\{a\}}$,
\begin{eqnarray*}
\label{eq:transprob2}
\begin{array}{lll}
q({\bf n},{\bf n}+e_1)          &=& \lambda_1,~~ n_1=1,\ldots,L-2,\\
q({\bf n},{\bf n}+e_l)          &=& \lambda_l,~~ l=2,\ldots,M,\\
q({\bf n},{\bf n}-e_1)          &=& 1/b_1,  ~~ 2 \leq n_1 \leq L-1, \\
q\big({\bf n},\{a\}\big)        &=& \alpha_1, ~~ 2 \leq n_1 \leq L-1, \\
q\big({\bf n},\{a\}\big)        &=& \alpha_1+1/b_1, ~~  n_1=1.
\end{array}
\end{eqnarray*}

We derive now $\mathbb{P}_L\big({\bf N}_1^e=(j_1,\ldots,j_M)~|~{\bf N}_1^b=(i_1,\ldots,i_M)\big)$.
We order the transient states lexicographically first according to $n_M$, then to $n_{M-1}$, $\ldots$, and
finally to $n_1$. Similarly to the
time-limited discipline, during a server visit to $Q_1$, the number of
jobs at $Q_j$, $j=2,\ldots,M$, may only increase. It then follows that
the transient generator of {\bf AMC}$_{T}$ has the same structure as the
transient generator of {\bf AMC}$_{A}$, i.e.~it is an upper-bidiagonal Toeplitz
matrix of upper-bidiagonal Toeplitz diagonal blocks. Therefore, by the same
arguments as for the time-limited discipline, we find that the joint moment of the
p.g.f.~of ${\bf N}^e_1$ and the event that the absorption is due to
timer expiration, denoted by $\{\mbox{timer}\}$, given ${\bf N}_1^b$, reads
\begin{eqnarray*}
\mathbb{E}_L\Big[{\bf z}^{{\bf N}^e_1}\cdot {\bf 1}_{\{\mbox{timer}\}}~|~{\bf N}_1^b=(i_1,\ldots,i_M) \Big]=\nonumber~~~~~~~~~~~~~~~~~~~~~~~~~~~~~~~~~~~~~~~~~&& \\
~~~~~~~~~~~~~~~~~~~~~~~~~ -\alpha_1 z_{2}^{i_{2}} \ldots z_M^{i_M} e_{i_1}\Big({\bf B}_{2}+(z_{2}\lambda_{2}+\ldots+z_M\lambda_M){\bf I}\Big)^{-1}g_{1}(z_1),&&
\label{eq:ztran0ET}
\end{eqnarray*}
where ${\bf B}_2$ is the sum of the generator matrix of an M/M/1/L-1 queue with arrival rate
$\lambda_1$ and service rate $1/b_1$ restricted to the states with the number of jobs strictly
positive, and of the matrix $-(\lambda_2+\ldots+\lambda_M+\alpha_1)${\bf I}, and where
\begin{eqnarray}
g_{1}(z_1)&:=&(z_1,\ldots,z_1^{L-1})^T.\nonumber
\end{eqnarray}
Let,
\begin{equation}
\label{eq:qetz}
{\bf Q_{T}(z)}:={\bf B}_{2}+(z_{2}\lambda_{2}+\ldots+z_M\lambda_M){\bf I}.
\end{equation}
The joint moment of the p.g.f.~of ${\bf N}^e_1$ and the event that the absorption is due to empty $Q_1$, denoted by $\{\mbox{$Q_1$ empty}\}$,
given ${\bf N}_1^b$, reads
\begin{eqnarray*}
\mathbb{E}_L\Big[{\bf z}^{{\bf N}^e_1}\cdot {\bf 1}_{\{\mbox{$Q_1$ empty}\}}~|~{\bf N}_1^b=(i_1,\ldots,i_M) \Big]=-1/b_1 z_{2}^{i_{2}} \ldots z_M^{i_M} e_{i_1}\big({\bf Q_{T}(z)}\big)^{-1}e_1,
\label{eq:ztran0ET1}
\end{eqnarray*}
Summing the latter two p.g.f.~gives the p.g.f.~of ${\bf N}^e_1$ given ${\bf N}_1^b$, which
reads
\begin{eqnarray}
\mathbb{E}_L\big[{\bf z}^{{\bf N}^e_1}~|{\bf N}_1^b=(i_1,\ldots,i_M)~\big]&=&-\alpha_1 z_{2}^{i_{2}} \ldots z_M^{i_M} e_{i_1} \times \nonumber \\
&&\big({\bf Q_{T}(z)}\big)^{-1}\Big(g_{1}(z_1)+\frac{1}{b_1\alpha_1}\cdot e_1\Big),
\label{eq:ztran1ET}
\end{eqnarray}

In the final part of this section, we find the inverse of ${\bf Q_{T}(z)}$ and let $L\to\infty$.
%\subsection{Inverse of ${\bf Q_{T}(z)}$}
%\label{}

We note that ${\bf Q_{T}(z)}={\bf T(z)}+\lambda_1 vv^T$, $v=(0,\ldots,0,1)^T$,
where ${\bf T_{T}(z)}$ is a (L-1)-by-(L-1) tridiagonal Toeplitz matrix with diagonal entries
equal to $\big(-\lambda_1-1/b_1-\alpha_1-\sum_{m=2}^{M}\lambda_m(1-z_m)\big)$, upper-diagonal
entries are equal to $\lambda_1$, and lower-diagonal entries $1/b_1$. We emphasize that the only
difference between ${\bf T_{A}(z)}$ of the autonomous-server discipline and ${\bf T_{T}(z)}$ is
that ${\bf T_{A}(z)}$ is an L-by-L matrix. Therefore, following the same approach as in
Section~\ref{sec:auto}, we find that
the $(i,j)$-entry of ${\bf Q_{T}(z)^{-1}}$, $i,j=1,\ldots,L-1$, gives
\begin{equation}
q(i,j)^*=t(i,j)^*_{T}-\lambda_1\frac{t(i,L-1)^*_{T}t(L-1,j)^*_{T}}{1+\lambda_1 t(L-1,L-1)^*_{T}},
\label{eq:invQETz1}
\end{equation}
where $t(i,j)^*_{T}$ is the (i,j)-entry of ${\bf T_{T}(z)^{-1}}$ that reads
\begin{equation}
t(i,j)^*_{T}=\left\{\begin{array}{l l}
-\frac{(r_{11}^i-r_{21}^i)(r_{11}^{L-j}-r_{21}^{L-j})}{\lambda_1(r_{11}-r_{21})(r_{11}^{L}-r_{21}^{L})}            &,~i\leq j\leq L-1\\
\frac{(r_{11}^{-j}-r_{21}^{-j})(r_{11}^{L}r_{21}^i-r_{21}^{L}r_{11}^i)}{\lambda_1(r_{11}-r_{21})(r_{11}^{L}-r_{21}^{L})} &,~j\leq i\leq L-1\\
\end{array} \right.
\label{eq:invTET}
\end{equation}
where $r_{11}$ and $r_{21}$ are the distinct roots of $P_1(r):=\lambda_1 r^2-s_1r+1/b_1$.
Inserting the values of $q(i,j)^*_{T}$ into~(\ref{eq:ztran1ET}) yields that
\begin{align}
\lefteqn{\mathbb{E}_L\big[{\bf z}^{{\bf N}^e_1}~|{\bf N}_1^b=(i_1,\ldots,i_M)~\big]=-\alpha_1 z_{2}^{i_{2}} \ldots z_M^{i_M} \times ~~~~~~~~~~~~~~~~~~~~~~~~~}\nonumber \\
&~~~~~~~~~~~~\bigg[\frac{1}{b_1\alpha_1}q(i_1,1)^*+\sum_{j=1}^{L-1}z_1^{j}\bigg(t(i_1,j)^*_{T}
-\lambda_1\frac{t(i_1,L-1)^*_{T}\cdot t(L-1,j)^*_{T}}{1+\lambda_1 t(L-1,L-1)^*_{T}} \bigg)\bigg].\nonumber\\
\label{eq:ztran2ET}
\end{align}
Some technical calculus shows that the following limit is equal to zero
\begin{eqnarray}
\lim_{L\to\infty}\frac{t(i_1,L-1)^*_{T}}{1+\lambda_1 t(L-1,L-1)^*_{T}}\sum_{j=1}^{L-1}z_1^{j}\cdot t(L-1,j)^*_{T}.\nonumber
\end{eqnarray}
Plugging the latter limit, $q(i_1,1)^*$, and $t(i_1,j)^*_{T}$ in (\ref{eq:ztran2ET}), we find that
\begin{eqnarray}
\mathbb{E}\big[{\bf z}^{{\bf N}^e_1}~|{\bf N}_1^b=(i_1,\ldots,i_M)~\big]=
 z_{2}^{i_{2}} \ldots z_M^{i_M}\bigg(r_{11}^{i_1}-\frac{\alpha_1z_1}{P_1(z_1)}(z_1^{i_1}-r_{11}^{i_1}) \bigg)
\label{eq:ztran2ET1}
\end{eqnarray}
Removing the condition of ${\bf N}_1^b=(i_1,\ldots,i_M)$ for $i_1=0,\ldots,$ $L-1$,
\begin{eqnarray}
%\mathbb{E}\big[{\bf z}^{{\bf N}^e_1}\big]&=&\Big(1+\frac{\alpha_1z_1}{P_1(z_1)}\Big)\mathbb{E}\big[({\bf z}_1^*)^{{\bf N}^b_1}\big]-\frac{\alpha_1z_1}{P_1(z_1)}\mathbb{E}\big[{\bf z}^{{\bf N}^b_1}\big],\nonumber \\
\gamma^{T}_1({\bf z})&=&\Big(1+\frac{\alpha_1z_1}{P_1(z_1)}\Big)\beta^{T}_1({\bf z}_1^*)-\frac{\alpha_1z_1}{P_1(z_1)}\beta^{T}_1({\bf z}),
\label{eq:ztransNET}
\end{eqnarray}
where ${\bf z}^*_1:=(r_{11},z_2,\ldots,z_M)$. From (\ref{eq:ztransNET}), we deduce that for a
server visit to $Q_i$, $i=1,\ldots,M$,
\begin{equation}
\gamma^{T}_i({\bf z})=\Big(1+\frac{\alpha_iz_i}{P_i(z_i)}\Big)\beta^{T}_i({\bf z^*})-\frac{\alpha_iz_i}{P_i(z_i)}\beta^{T}_i({\bf z}),
\end{equation}
where ${\bf z}^*_i=(z_1,\ldots,z_{i-1},r_{1i},z_{i+1},\ldots,z_M)$, $|r_{1i}|<1$, and where $P_i(z_i)$, $s_i$, and $r_{1i}$ are in
(\ref{eq:pizi}), (\ref{eq:si}), and (\ref{eq:r1i}) respectively.
%\begin{eqnarray*}
%P_i(z_i) &=& \lambda_i z_i^2-s_iz_i+1/b_i, \\
%s_i      &=& \lambda_i+1/b_i+\alpha_i+\sum_{m=1,m\neq i}^{M}\lambda_m(1-z_m),\\
%r_{1i}   &=& \frac{s_i-\sqrt{(s_i)^2-4\lambda_i/b_i}}{2\lambda_i},
%\end{eqnarray*}

Finally, introducing the switch-over times from $Q_{i-1}$ to $Q_{i}$, we obtain
\begin{eqnarray}
\gamma^{T}_i({\bf z})&=&\Big(1+\frac{\alpha_iz_i}{P_i(z_i)}\Big)\gamma^{T}_{i-1}({\bf z}^*_i)C^{i-1}({\bf z}^*_i)
                     -\frac{\alpha_iz_i}{P_i(z_i)}\gamma^{T}_{i-1}({\bf z})C^{i-1}({\bf z}).
\end{eqnarray}

\section{k-Limited Discipline}
\label{sec:k-limit}
In this section, we analyze the $k$-limited discipline.
According to this discipline the server continues working at a queue until
either a predefined number of $k$ jobs is served or the queue becomes
empty, whichever occurs first. Similarly to the previous disciplines, the objective
is to relate the joint queue-length probabilities at the
beginning and end of a server visit to $Q_1$, referred to as
$\beta_1^{k}({\bf z})$ and $\gamma_1^{k}({\bf z})$.

By analogy with the time-limited discipline, we will first assume that $Q_1$
has a limited queue of $L-1$ jobs, second there are ${\bf N}_1^b:=(i_1,...,i_M)$
jobs in $(Q_1,\ldots,$ $Q_M)$, with $i_1\geq 1$, at the beginning time of a server
visit to $Q_1$, and third there are ${\bf N}_1^b:=(j_1,...,j_M)$ jobs in $(Q_1,\ldots,Q_M)$
at the end time of a server visit to $Q_1$. Note that if $Q_1$ is empty at the beginning
of a server visit, $i_1=0$, the server will leave immediately, i.e.,
$\mathbb{P}\big({\bf N}_1^e={\bf N}_1^b\big)=1$. For
this reason, we will exclude the latter obvious case from the analysis in the following,
however, we will include it when we will uncondition on ${\bf N}_1^b$.

Let ${\bf N}(t):=(N_1(t),\ldots,N_M(t),D(t))$ denote the $M+1$-dimensional,
continuous-time Markov chain with discrete state-space
$\xi_{k}=\{1,\ldots,L-1\}\times\{0,1,\ldots\}^{M-1}\times\{0,1,\ldots\}\cup\{a\}$, where
$N_j(t)$ represents the number of jobs in $Q_j$ at time $t$ during a server visit to $Q_1$,
and $D(t)$ is the total number of departures from $Q_1$ until $t$. State $\{a\}$ is absorbing.
This absorbing Markov chain is denoted by ${\bf AMC}_{k}$. The absorption of ${\bf AMC}_{k}$
occurs when the server leaves $Q_1$ which happens with rate $1/b_1$ from all
transient states with $D(t)=k-1$ or $N_1(t)=1$.

We set ${\bf N}(0)=({\bf N}_1^b,0)$. The probability that the transition
to the absorption state occurs from one of the states
$\{(j_1,\ldots,j_M)\}$, $j_1\geq 2$, equals
$\mathbb{P}_L\big({\bf N}_1^e=(j_1-1,\ldots,j_M)~|~{\bf N}_1^b\big)$ and the absorption is
eventually due to $k$ departures from $Q_1$ with rate $1/b_1$. If the
absorption is due to $Q_1$ becoming empty, $\mathbb{P}_L\big({\bf N}_1^e=(0,j_2\ldots,j_M)~|~{\bf N}_1^b\big)$
equals the probability that the transition to absorption is with rate $1/b_1$ and it
occurs from state $\{(1,j_2,\ldots,j_M)\}$. Note that it is possible that the k-th
departure at $Q_1$ leaves behind an empty queue. In our analysis we will consider
this event as a transition to absorption that is due to $Q_1$ becoming empty.
The non-zero transition rates of ${\bf AMC}_{k}$ can be written for all
${\bf n}=(n_1,\ldots,n_M,j)$ $\in\xi_{k}-\{a\}$,

\begin{eqnarray*}
\label{eq:transprob3}
\begin{array}{lll}
q({\bf n},{\bf n}+e_1)        &=& \lambda_1,~~~~n_1=1,\ldots,L-2,\\
q({\bf n},{\bf n}+e_l)        &=& \lambda_l,~~~~l=2,\ldots,M,\\
q({\bf n},{\bf n}-e_1+e_{M+1})&=& 1/b_1,~n_1=2,\ldots ,L-1,\\
&&~~~~~~~~j=0,\ldots ,k-2,\\
q({\bf n},\{a\})              &=& 1/b_1,     ~n_1=1~\mbox{or}~j= k-1.
\end{array}
\end{eqnarray*}

We derive now $\mathbb{P}_L\big({\bf N}_1^e=(j_1,\ldots,j_M)~|~{\bf
  N}_1^b=(i_1,\ldots,i_M)\big)$. We order the transient states of
${\bf AMC}_{k}$ lexicographically according to $n_M$, $n_{M-1}$,
$\ldots$, $n_2$, then to $j$, and finally according to $n_1$. During a
server visit to $Q_1$, the number of jobs at $Q_j$, $j=2,\ldots,M$, may only increase. Therefore,
similarly to the automous-server and time-limited discipline, we deduce that the
joint moment of the p.g.f.~of ${\bf N}^e_1$ and the event that the absorption
is due $k$ to departures, denoted by $\{\mbox{k~dep.}\}$, given ${\bf N}_1^b$, reads
\begin{eqnarray}
\mathbb{E}_L\Big[{\bf z}^{{\bf N}^e_1} {\bf 1}_{\{\mbox{k~dep.}\}}~|~{\bf N}_1^b=(i_1,\ldots,i_M) \Big]=-1/b_1 z_{2}^{i_{2}} \ldots z_M^{i_M}e_{i_1}\otimes e_{1} \times ~~~~~~~~~\nonumber&& \\
 \Big({\bf C}_{2}+(z_{2}\lambda_{2}+\ldots+z_M\lambda_M){\bf I}\Big)^{-1}h(z_1),&&
\label{eq:ztran0k}
\end{eqnarray}
where $e_1$ is a k-dimensional row vector of zero entries except the first that is
one, $({\bf C}_{2}+(z_{2}\lambda_{2}+\ldots+z_M\lambda_M){\bf I})$ is a k-by-k
upper-bidiagonal block matrix of upper diagonal blocks equal to {\bf U}, where {\bf U}
is an (L-1)-by-(L-1) lower-diagonal matrix whose entries equal to $1/b_1$,
and of diagonal blocks equal to {\bf D}, where {\bf D} is the sum of the generator matrix of a
M/M/1/L-1 queue with arrival rate $\lambda_1$ and service rate $0$ restricted to strictly
positive states, and of the matrix $-(\lambda_2(1-z_2)+\ldots+\lambda_M(1-z_M)+1/b_1)${\bf I}, and
\begin{eqnarray}
h(z_1)&:=&q(z_1)\otimes e_{k},\\
q(z_1)&:=&(0,z_1,\ldots,z_1^{L-2})^T,\label{eq:oz1}
\end{eqnarray}
and where $e_{k}$ is a k-dimensional column vector of zero entries except the k-th that is one.
Plugging the inverse of $\Big({\bf C}_{2}+(z_{2}\lambda_{2}+\ldots+z_M\lambda_M){\bf I}\Big)$
into (\ref{eq:ztran0k}) gives that
\begin{eqnarray}
\mathbb{E}_L\Big[{\bf z}^{{\bf N}^e_1}\cdot {\bf 1}_{\{\mbox{k~dep.}\}}~|~{\bf N}_1^b=(i_1,\ldots,i_M) \Big]&=& -1/b_1 z_{2}^{i_{2}} \ldots z_M^{i_M}e_{i_1} \times \nonumber \\
&&\big(-{\bf D}^{-1}{\bf U}\big)^{k-1}{\bf D}^{-1}q(z_1),
\label{eq:ztran1k}
\end{eqnarray}
The joint moment of the p.g.f.~of ${\bf N}^e_1$ and the event that the absorption is due to empty $Q_1$, denoted by
$\{\mbox{$Q_1$ emp}\}$, given ${\bf N}_1^b$, reads
\begin{eqnarray}
\lefteqn{\mathbb{E}_L\Big[{\bf z}^{{\bf N}^e_1} {\bf 1}_{\{\mbox{$Q_1$ emp.}\}}~|~{\bf N}_1^b=(i_1,\ldots,i_M) \Big]~~~~~~~~}\nonumber\\
%\end{eqnarray}
%\begin{eqnarray}
&=&-1/b_1 z_{2}^{i_{2}} \ldots z_M^{i_M} e_{i_1}\otimes e_{1}
 \Big({\bf C}_{2}+(z_{2}\lambda_{2}+\ldots+z_M\lambda_M){\bf I}\Big)^{-1}e_1\otimes e \nonumber\\
&=&-1/b_1 z_{2}^{i_{2}} \ldots z_M^{i_M}e_{i_1}\big({\bf I}-\big(-{\bf D}^{-1}{\bf U}\big)^{k}\big)\big({\bf D}+{\bf U}\big)^{-1}e_1.
\label{eq:ztran2k}
\end{eqnarray}
Summing the latter two p.g.f.~gives $\mathbb{E}_L\big[{\bf z}^{{\bf N}^e_1}~|~{\bf N}_1^b\big]$. It remains to find
first $e_{i_1}\big(-{\bf D}^{-1}{\bf U}\big)^{k}$, second $\big({\bf D}^{-1}\big)q(z_1)$ and $\big({\bf D}+{\bf U}\big)^{-1}e_1$,
so that finally we will take the limit for $L\to\infty$ of $\mathbb{E}_L\big[{\bf z}^{{\bf N}^e_1}~|~{\bf N}_1^b\big]$ .

\subsection{$e_{i_1}\big(-{\bf D}^{-1}{\bf U}\big)^{k}$}
\label{sec:Mk}

\begin{figure}[hbtp]
\vspace{0pt}
  \centerline{\hbox{
    \includegraphics[width=3.2in]{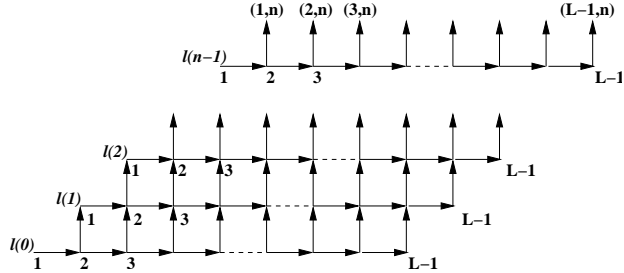}
  }}

\vspace{0cm}
\caption{Directed graph for the computation of $c^n(i,j)$.}
  \label{fig:dgraph}
\end{figure}

The matrix {\bf D} is an (L-1)-by-(L-1) upper-bidiagonal matrix with upper-diagonal entries equal to
$\lambda_1$ and diagonal equal to $-\lambda_1(x,\ldots,x,x_0)$, where $x:=(\lambda_1+\lambda_2(1-z_2)+\ldots+\lambda_M(1-z_M)+1/b_1)/\lambda_1$ and
$x_0:=(\lambda_2(1-z_2)+\ldots+\lambda_M(1-z_M)+1/b_1)/\lambda_1$. Thus, it is easy to show that
$-{\bf D}^{-1}{\bf U}=(xb_1\lambda_1)^{-1}{\bf L}$, where
\begin{eqnarray*}
{\bf L}=\left(
\begin{array}{cccccccc}
x^{-1} & x^{-2} & x^{-3} & \cdots & x^{-L+3} & x^{-L+3}x_0^{-1} & 0 \\
1      & x^{-1} & x^{-2} & \cdots & x^{-L+4} & x^{-L+4}x_0^{-1} & 0 \\
0      & 1      & x^{-1} & \cdots & x^{-L+5} & x^{-L+5}x_0^{-1} & 0 \\
\vdots & \ddots & \ddots & \ddots & \vdots   & \vdots           & \vdots \\
0      & \cdots & \cdots & 0      & 1        & x_0^{-1}         & 0 \\
0      & \cdots & \cdots & \cdots & 0 & x\cdot x_0^{-1} & 0
\end{array}
\right).
\end{eqnarray*}
For $n\geq 1$, note that the (i,j)-entry of ${\bf L}^n$, can be written
as $c^n(i,j)x^{-n+i-j}$, $j=1,\ldots,L-3$. We do not consider the (i,j)-entry
of ${\bf L}^n$ with $j\geq L-2$ since these entries will tend to zero when
we will take the limit for $L\to\infty$.
The coefficients $c^n(i,j)$ are strictly positive integers for
$1\leq i\leq n$ and $1 \leq j \leq L-2$, and $n+1 \leq i \leq L-1$ and
$i-n \leq j \leq L-2$, and zero otherwise. Moreover, the
sequence $\{c^n(i,j)\}$ satisfies the following recurrent equation for $n\geq 2$,
$1\leq i \leq L-1$ and $1\leq j\leq L-2$,
\begin{eqnarray}
c^n(i,j)=c^n(i,j-1)+c^{n-1}(i,j+1)=\sum_{l=1}^{j+1} c^{n-1}(i,l),
\end{eqnarray}
where
\begin{eqnarray}
c^1(i,j)=\left\{
\begin{array}{lll}
1, &i=1, 1\leq j \leq L-2, &\\
1, &2 \leq i \leq L-1, i-1 \leq j \leq L-2,&\\
0, &\mbox{otherwise}.&
\end{array}
\right.
\end{eqnarray}
The coefficient $c^n(i,j)$ can be interpreted as the number of paths in the directed graph in
Figure~\ref{fig:dgraph}. Especially, $c^n(i,j)$ equals the number of paths from state $i$ in level
$l(0)$ to state $j$ in level $l(n)$. Thus by an induction argument, it can be shown that
$c^n(i,j)$ has the following solution for $2\leq n<L-2$. That is, for $j=1,\ldots,L-2$ and $n\ll L$,
\begin{equation}
c^n(1,j)=\left(\begin{array}{c}2n+j-2\\n-1\end{array}\right)-\left(\begin{array}{c}2n+j-2\\n+j\end{array}\right),
\end{equation}
for $i=2,\ldots,n-1$ and $j=1,\ldots,L-2$,
\begin{equation}
c^n(i,j)=\left(\begin{array}{c}2n+j-i-1\\n-1\end{array}\right)-\left(\begin{array}{c}2n+j-i-1\\n+j\end{array}\right),
\end{equation}
for $i=n$ and $j=1,\ldots,L-2$,
\begin{eqnarray}
c^n(n,j)&=&\left(\begin{array}{c}n+j-1\\n-1\end{array}\right),
\end{eqnarray}
for $i=1,\ldots,L-n-1$ and $j=i,\ldots,L-2$,
\begin{eqnarray}
c^n(i+n,j)&=&\left(\begin{array}{c}n+j-i-1\\n-1\end{array}\right),
\end{eqnarray}
and $c^n(i,j)$ equals zero for $i=n+2,\ldots,L-1$ and $j=1,\ldots,i-n$.

Finally, we conclude that $e_{i_1}\big(-{\bf D}^{-1}{\bf U}\big)^{k}$ is a row vector of size $L-1$ that is
equal to $(xb_1\lambda_1)^{-k}{\bf L}^k$ with $j$-th element equal to
\begin{equation}
\frac{c^k(i_1,j)}{(\lambda_1 b_1)^k}x^{-2k+i_1-j},
\label{eq:entries}
\end{equation}
for $j=1,\ldots,L-3$. Note that since $|x|<1$, the limit of (\ref{eq:entries}) tends zero for $L\to\infty$.

\subsection{$\big({\bf D}^{-1}\big)q(z_1)$ and $\big({\bf D}+{\bf U}\big)^{-1}e_1$}

The matrix {\bf D} is a (L-1)-by-(L-1) upper-bidiagonal matrix with upper-diagonal $(\lambda_1,\ldots,\lambda_1)$ and
diagonal $-\lambda_1(x,\ldots,x,x_0)$, where $x$ and $x_0$ are defined in Section~\ref{sec:Mk}. Thus,
\begin{eqnarray*}
{\bf D}^{-1}=\lambda_1^{-1}\left(
\begin{array}{cccccccc}
x^{-1} & x^{-2} & x^{-2} & \cdots & x^{-L+2} & x^{-L+2}x_0^{-1}\\
0      & x^{-1} & x^{-2} & \cdots & x^{-L+3} & x^{-L+3}x_0^{-1}\\
0      & 0      & x^{-1} & \cdots & x^{-L+4} & x^{-L+4}x_0^{-1}\\
\vdots & \ddots & \ddots & \ddots & \vdots   & \vdots          \\
0      & \cdots & \cdots & 0      & x^{-1}   & x^{-1} x_0^{-1} \\
0      & \cdots & \cdots & \cdots & 0        &        x_0^{-1}
\end{array}
\right).
\end{eqnarray*}
Using (\ref{eq:oz1}), we find that ${\bf D}^{-1}q(z_1)$ is an (L-1)-dimensional column vector of
$i$-th element, denoted as $d(i)$, equal to
\begin{eqnarray}
d(1)&=&-\frac{1}{\lambda_1}\Big(z_1 x^{-1}\frac{1-(z_1x^{-1})^{L-3}}{x-z_1}+x^{-L+2}x_0^{-1}z_1^{L-2}\Big),\nonumber\\
d(i)&=&-\frac{1}{\lambda_1}\Big(z_1^{i-1} \frac{1-(z_1x^{-1})^{L-1-i}}{x-z_1}+x^{-L+i+1}x_0^{-1}z_1^{L-2}\Big),\nonumber\\
\label{eq:dz1}
\end{eqnarray}
for $i=2,\ldots,L-1$. Note that $|z_1/x|<1$ which gives that
\begin{eqnarray}
\lim_{L\to\infty}d(1)&=&-\frac{1}{\lambda_1}\Big(\frac{z_1 x^{-1}}{x-z_1}\Big),\\
\lim_{L\to\infty}d(i)&=&-\frac{1}{\lambda_1}\Big(\frac{z_1^{i-1}}{x-z_1}\Big),
\label{eq:dz1limit}
\end{eqnarray}
for all $i<\infty$.

Now we compute $\big({\bf D}+{\bf U}\big)^{-1}e_1$. Recall that
$\big({\bf D}+{\bf U}\big)$ is an (L-1)-by-(L-1) tridiagonal matrix
with upper-diagonal entries equal $\lambda_1$, diagonal $-\lambda_1(x,\ldots,x,x_0)$
and lower-diagonal entries $1/b_1$. Therefore, $\big({\bf D}+{\bf U}\big)$ is equal to
the matrix ${\bf Q_{T}(z)}$ in (\ref{eq:qetz}) with $\alpha_1=0$. We note that the
inverse of ${\bf Q_{T}(z)}$ was computed in (\ref{eq:invQETz1}), thus using these results we find that $\big({\bf
  D}+{\bf U}\big)^{-1}e_1$ is a column vector equal to $\big(p(1),\ldots,p(L-1)\big)^T$ with the $i$-th entry that is given by
\begin{eqnarray}
p(i)&:=&t(i,1)_T^*-\lambda_1\frac{t(i,L-1)_T^* t(L-1,1)_T^*}{1+\lambda_1t(L-1,L-1)_T^*},\nonumber\\
   % &=&-b_1\frac{\big(\frac{r_{11}}{r_{21}}\big)^Lr_{21}^i-r_{11}^i}{\big(\frac{r_{11}}{r_{21}}\big)^L - 1 }+
%\lambda_1 b_1^2\frac{(r_{11}-r_{21})(r_{11}^i-r_{21}^i)}{ (\frac{1}{r_{21}^L}-\frac{1}{r_{11}^L}) (r_{11}^L-r_{11}^{L-1}-r_{21}^L+r_{21}^{L-1}) }\\
    &=&-b_1\frac{y_{11}^Ly_{21}^{i-L}-y_{11}^i}{y_{11}^Ly_{21}^{-L} - 1 }+\lambda_1 b_1^2\frac{y_{11}-y_{21}}{y_{21}^{-L}-y_{11}^{-L}}\times\nonumber\\
&&\frac{y_{11}^i-y_{21}^i}{ y_{11}^L-y_{11}^{L-1}-y_{21}^L+y_{21}^{L-1}}.
\label{eq:invDplusU}
\end{eqnarray}
where $y_{11}$ and $y_{21}$ are the distinct roots of
\begin{equation}
\lambda_1 y^2-s^*_1y+1/b_1,
\end{equation}
where $s^*_1:=\lambda_1+1/b_1+\sum_{m=2}^{M}\lambda_m(1-z_m)$. Note
that in this case $|y_{11}|\leq 1<|y_{21}|$, so that we may find that,
\begin{eqnarray}
\lim_{L\to\infty}p(i) &=&-b_1y_{11}^i,
\label{eq:invDplusULimit}
\end{eqnarray}
for all $i<\infty$.

\subsection{Limit of $\mathbb{E}_L\big[{\bf z}^{{\bf N}^e_1}~|~{\bf N}_1^b\big]$ for $L\to\infty$}

Plugging $e_{i_1}\big({\bf D}^{-1}{\bf U}\big)^{k}$, $\big({\bf D}^{-1}\big)q(z_1)$ and
$\big({\bf D}+{\bf U}\big)^{-1}e_1$ into $\mathbb{E}_L\big[{\bf z}^{{\bf N}^e_1}~|~{\bf N}_1^b\big]$
and taking the limit for $L\to\infty$ gives that
\begin{eqnarray}
\mathbb{E}\big[{\bf z}^{{\bf N}^e_1}~|~{\bf N}_1^b\big]&=&-1/b_1 z_{2}^{i_{2}}\ldots z_M^{i_M} S,\nonumber
\end{eqnarray}
where,
\begin{eqnarray}
S&:=&\frac{b_1c^{k-1}(i_1,1)}{(\lambda_1b_1)^{k}}x^{-2k+i_1}-b_1y_{11}^{i_1} -\frac{b_1x^{-2k+i_1+1}}{(\lambda_1b_1)^{k}(x-z_1)}\sum_{j=1}^{\infty}c^{k-1}(i_1,j)\Big(\frac{z_1}{x}\Big)^{j-1}\nonumber\\
&&+\frac{b_1x^{-2k+i_1}}{(\lambda_1 b_1)^{k}}\sum_{j=1}^{\infty}c^{k}(i_1,j)\Big(\frac{y_{11}}{x}\Big)^j,
\label{eq:S}
\end{eqnarray}
for $k\geq 2$.

Due to the complexity of the analysis for an arbitrary $k$, we will restrict ourselves to the
$1$-limited and $2$-limited disciplines.
\newline\\
{\bf $\bf 1$-limited}. First take the limits of $d(i)$ and $p(i)$ in (\ref{eq:dz1limit}) and (\ref{eq:invDplusULimit}), then
plugging $k=1$ into (\ref{eq:ztran1k}) and (\ref{eq:ztran2k}) gives that
\begin{eqnarray}
\mathbb{E}\big[{\bf z}^{{\bf N}^e_1}~|~{\bf N}_1^b\big]&=&%z_2^{i_2}\ldots z_{M}^{i_M}
%\Big(\frac{z_1x^{-1}}{\lambda_1 b_1(x-z_1)}+y_{11}-\frac{y_{11}x^{-1}}{\lambda_1 b_1(x-y_{11})}\Big),\nonumber\\
\frac{z_2^{i_2}\ldots z_{M}^{i_M}}{\lambda_1 b_1 x}\bigg(\frac{z_1}{x-z_1}+1\bigg),
\end{eqnarray}
for $i_1=1$, and
\begin{eqnarray}
\mathbb{E}\big[{\bf z}^{{\bf N}^e_1}~|~{\bf N}_1^b\big]&=&%z_2^{i_2}\ldots z_{M}^{i_M}
%\Big(\frac{z_1^{i_1-1}}{\lambda_1 b_1(x-z_1)}+y_{11}^{i_1}-\frac{y_{11}^{i_1-1}}{\lambda_1 b_1(x-y_{11})}\Big),\nonumber\\
z_2^{i_2}\ldots z_{M}^{i_M}\Big(\frac{z_1^{i_1-1}}{\lambda_1 b_1(x-z_1)}\Big),
\end{eqnarray}
for $i_1=2,3,\ldots$~. Unconditioning on ${\bf N}_1^b=({\bf N}_{11}^b,\ldots,{\bf N}_{M1}^b)$, we find that
\begin{eqnarray}
\mathbb{E}\big[{\bf z}^{{\bf N}^e_1}\big]=\frac{1/b_1z_1^{-1}}{1/b_1+\lambda_1(1-z_1)+\lambda_2(1-z_2)}\mathbb{E}\big[{\bf z}^{{\bf N}^b_1}\big]+\nonumber\\
\bigg(1-\frac{1/b_1z_1^{-1}}{1/b_1+\lambda_1(1-z_1)+\lambda_2(1-z_2)}\bigg)\mathbb{E}\big[{\bf z}^{{\bf N}^b_1}\big]\Big|_{z_1=0}.
\end{eqnarray}
\newline\\
{\bf $\bf 2$-limited}.  Plugging $k=2$ in (\ref{eq:S}) gives that p.g.f.~of ${\bf N}^e_1$ then gives
\begin{eqnarray}
\mathbb{E}\big[{\bf z}^{{\bf N}^e_1}~|~{\bf N}_1^b\big]%&=&z_2^{i_2}\ldots z_{M}^{i_M}
%\Big(\frac{z_1x^{-3}(2x-z_1)}{\lambda_1^2 b_1^2(x-z_1)^2}+y_{11}-\frac{y_{11}^3}{x^{3}}(x+y_{21})\Big),\nonumber\\
&=&\frac{z_2^{i_2}\ldots z_{M}^{i_M}}{\lambda_1b_1 x}
\Big(\frac{1}{\lambda_1b_1(x-z_1)^2}+1\Big),
\end{eqnarray}
for $i_1=1$, and
\begin{equation}
\mathbb{E}\big[{\bf z}^{{\bf N}^e_1}~|~{\bf N}_1^b\big]=z_2^{i_2}\ldots z_{M}^{i_M}
\Big(\frac{z_1^{i-2}}{\lambda_1^2b_1^2(x-z_1)^2}\Big),
\end{equation}
for $i=2,3,\ldots$~. Unconditioning on ${\bf N}_1^b=({\bf N}_{11}^b,\ldots,{\bf N}_{M1}^b)$, we find that
\begin{eqnarray}
\mathbb{E}\big[{\bf z}^{{\bf N}^e_1}\big]&=&\frac{z_1^{-2}}{\lambda_1^2b_1^2(x-z_1)^2}\mathbb{E}\big[{\bf z}^{{\bf N}^b_1}\big]+
\bigg(1-\frac{z_1^{-2}}{\lambda_1^2b_1^2(x-z_1)^2}\bigg) \mathbb{E}\big[{\bf z}^{{\bf N}^b_1}\big]\Big|_{z_1=0} \nonumber\\
&&+\bigg(\frac{z_1^{-1}}{\lambda_1 b_1 x}-\frac{z_1^{-2}}{\lambda_1^2 b_1^2 x (x-z_1)}\bigg)\mathbb{E}\big[{\bf 1}_{\{N_{11}^b=1\}}{\bf z}^{{\bf N}^b_1}\big].
%+\bigg(-\frac{1}{\lambda_1^2b_1^2x^2}+y_{11}^2-\frac{y_{11}^3}{x^{2}}(x+y_{21})\bigg)\mathbb{E}\big[{\bf 1}_{\{N_{11}^b=2\}}{\bf z}^{{\bf N}^b_1}\big].&&\nonumber\\
\end{eqnarray}

\begin{remark}
The results for 1-limited and 2-limited can also be obtained more directly by explicitly conditioning
on the number of jobs at the beginning of a server visit to a queue and keeping track how the queue-length
evolves. However, our analysis above shows that our tool can also applied to the k-limited discipline for
$k\geq 3$.
\end{remark}

\begin{remark}
{\bf Exhaustive discipline:}. The k-limited discipline for $k\to\infty$ is equivalent to the exhaustive discipline. Since
$\big(-{\bf D}^{-1}{\bf U}\big)$ is a sub-stochastic matrix with the sum of its row entries strictly
smaller than one, the limit $\big(-{\bf D}^{-1}{\bf U}\big)^k\to 0$ for $k\to\infty$. Therefore, taking the limit in
(\ref{eq:ztran1k}) and (\ref{eq:ztran2k}) for $k\to\infty$ and summing these limits give that
\begin{eqnarray}
\mathbb{E}_{L}\big[{\bf z}^{{\bf N}^e_1}~|~{\bf N}_1^b\big]=-1/b_1 z_{2}^{i_{2}}\ldots z_M^{i_M}e_{i_1}\big({\bf D}+{\bf U}\big)^{-1}e_1.
\end{eqnarray}
The limit of $\mathbb{E}_{L}\big[{\bf z}^{{\bf N}^e_1}~|~{\bf N}_1^b\big]$ for $L\to\infty$ then reads
\begin{eqnarray}
\mathbb{E}\big[{\bf z}^{{\bf N}^e_1}~|~{\bf N}_1^b\big]&=&y_{11}^{i_1}z_{2}^{i_{2}}\ldots z_M^{i_M}.
\end{eqnarray}
Finally, the unconditioning on ${\bf N}_1^b$ gives that
\begin{eqnarray}
\mathbb{E}\big[{\bf z}^{{\bf N}^e_1}\big]&=&\mathbb{E}\big[({\bf z}_1^*)^{{\bf N}^e_1}\big], \nonumber\\
\gamma^{E}_1({\bf z})&=&\beta_1^{E}({\bf z}_1^*),
\end{eqnarray}
where ${\bf z}_1^*=(y_{11},z_2,\ldots,z_M)$. Considering a server visit to $Q_i$, an equivalent
relation can be derived for $\gamma^{E}_i({\bf z})$ and $\beta_1^{E}({\bf z}_i^*)$ as follows
\begin{eqnarray}
\gamma^{E}_i({\bf z})&=&\beta_i^{E}({\bf z}_i^*).\nonumber
\end{eqnarray}
Now including $C^{i-1}$, the switch-over time from $Q_{i-1}$ and $Q_i$, it is easy to find
that
\begin{eqnarray}
\gamma^{E}_i({\bf z})&=&\gamma_{i-1}^{E}({\bf z}_1^*)C^{i-1}({\bf z}^*_i).\label{eq:exhau}
\end{eqnarray}
where ${\bf z}^*_i:=(z_1,\ldots,z_{i-1},y_{1i},z_{i+1},\ldots,z_M)$ and $y_{1i}$ is the
root of
\begin{equation}
\lambda_i y^2-s^*_iy+1/b_i,
\end{equation}
with $|y_{1i}|\leq 1$ and where $s^*_i=\lambda_i+1/b_i+\sum_{m=1,m\neq i}^{M}\lambda_m(1-z_m)$.
Eq.~(\ref{eq:exhau}) is equivalent to the well-known relation of exhaustive discipline in
(see, e.g., \cite[Eq.~(24)]{Eisenberg}).
\end{remark}

\section{Iterative scheme}
\label{sec:ite}

In this section, we will explain how to obtain the joint queue-length distribution using an iterative scheme.
First, let see how to compute $\gamma_i({\bf z})$ as function $\gamma_{i-1}({\bf z})$, where
${\bf z}=(z_1,\ldots,z_M)$.

Note that $\gamma_i({\bf z})$ is a function of $\gamma_{i-1}({\bf z})$ and
$\gamma_{i-1}({\bf z}_{i}^*)$ where ${\bf z}_i^*=(z_1,\ldots$ $,z_{i-1},a,z_{i+1}\ldots,z_M)$
with $|z_i|=1$, $i=1,\ldots,M$ and $|a|\leq 1$, which is a function of $z_l$ for all $l=1,\ldots,M$ and
$l\neq i$. Since
$\gamma_{i-1}({\bf z})$ is a joint p.g.f., the function  $\gamma_{i-1}({\bf z})$ is analytic
in $z_i$ for all $z_1,\ldots,z_{i-1},$ $z_{i+1},\ldots,z_M$. Hence, we can write
\[
\gamma_i({\bf z})=\sum_{n=0}^\infty g_{in}(z_1,\ldots,z_{i-1},z_{i+1}\ldots,z_M)z_i^n,\mbox{\hskip 1cm}|a|\le1,
\]
where $g_{in}(.)$ is again an analytic function. From complex function theory, it is well known that
\[
\gamma_{i}({\bf z}_{i}^*)=\frac 1{2\pi\imi}\oint_{C} \frac{\gamma_{i}({\bf z})}{z_i-a}dz_i,\mbox{\hskip 1cm for }|a|\le1,
\]
where $C$ is the unit circle and $\imi^2=-1$, and furthermore
\[
g_{in}(z_1,\ldots,z_{i-1},z_{i+1}\ldots,z_M)=\frac 1{2\pi\imi}\oint_{C} \frac{\gamma_{i}({\bf z})}{z_i^{n+1}}dz_i,\mbox{\hskip 1cm}
\]
where $n=0,1,\ldots~$. These formulas show that we only need to know
the joint p.g.f.\ $\gamma_{i-1}({\bf z})$ for all $\vecz$ with $|z_i|=1$,
to be able to compute  $\gamma_{i}({\bf z})$.% For this reason, one can
%use the Discrete Fourier Transform (d.f.t.) instead of the p.g.f., thus
%replacing $z_i$ by $w_i=e^{-2\pi{\bf i}k_i/J_i}$, where ${\bf i}^2=-1$, $k_i$
%refers to the i-th discrete point used for $Q_i$ in the d.f.t., and
%$J_i$ is the total number of latter points.

When there is an incurred switch-over time from queue $i-1$ to $i$
the p.g.f.~of the joint queue-length at the end of the n-th server
visit to $Q_i$, denoted by $\gamma^{n}_{i}({\bf z})$, can be computed
as a function of $\gamma^{n}_{i-1}({\bf z})$. The main step is to iterate
over all queues in order to express $\gamma^{n+1}_i({\bf z})$ as a
function of $\gamma^{n}_i({\bf z})$. Assuming that the
system is in steady-state these two latter quantities should be equal.
Thus, starting with an empty system at the first service visit to $Q_i$
and repeating the latter main step one can compute
$\gamma^2_i({\bf z})$, $\gamma^3_i({\bf z})$, and so on.  This iteration
is stopped when $\gamma^n_i({\bf z})$ converges.

\section{Tandem model}
\label{sec:tand}

We know that our tool can be applied also for Jackson-like queueing
networks with a single server that can serve only one queue at a
time. To show this, we will consider the example of a tandem model of
$M$ queues in series. $Q_1$ has Poisson arrivals. The service
requirement at $Q_i$ is distributed exponentially with mean $1/b_i$.
In the model there is only one server serving the queues according to
some schedule. The service discipline is either the autonomous-server
or the time-limited.  Observe that this tandem model is equivalent to
polling system with the property that only $Q_1$ has a Poisson
arrivals, the departures from $Q_i$ will join $Q_{i+1}$,
$i=1,\ldots,M-1$, and that departures
from $Q_{M}$ leaves the system.\\
\newline {\bf Autonomous-server}: according to this discipline the
server continues the service of a queue until certain exponentially
distributed time of rate $\alpha_1$ will elapse. Consider a server
visit to $Q_1$ following the same approach in Section~\ref{sec:auto}
we find that the solution is similar to (\ref{eq:ztran1}) and the matrix involved
has the same structure as ${\bf A}_2$. For this reason, we find that
\begin{equation}
\mathbb{E}\big[{\bf z}^{{\bf N}^e_1}\big]=\frac{\alpha_1(z_2-z_1)}{P_1(z_1)}\bigg(\frac{r_{11}\mathbb{E}\big[({\bf z}^*_1)^{{\bf N}^b_1}\big]}{z_2-r_{11}}-\frac{z_1\mathbb{E}[{\bf z}^{{\bf N}^b_1}]}{z_2-z_1}\bigg),
\end{equation}
where ${\bf z}^*_1:=(r_{11},z_2,\ldots,z_M)$ and  $r_{11}$ is the root of $P_1(r)=\lambda_1r^2-(\lambda_1+1/b_1+\alpha_1)r+z_2/b_1$ such that $|r_{11}|<1$. To relate $\mathbb{E}\big[{\bf z}^{{\bf N}^e_i}\big]$ to $\mathbb{E}\big[{\bf z}^{{\bf N}^b_i}\big]$ for a server visit to $Q_i$, $i>1$, we find that
\begin{equation}
\mathbb{E}\big[{\bf z}^{{\bf N}^e_i}\big]=\frac{\alpha_i(z_{i+1}-z_i)}{P_i(z_i)}\bigg(\frac{r_{1i}\mathbb{E}\big[({\bf z}_i^*)^{{\bf N}^b_i}\big]}{z_{i+1}-r_{1i}}-\frac{z_i\mathbb{E}[{\bf z}^{{\bf N}^b_i}]}{z_{i+1}-z_i}\bigg),
\end{equation}
where ${\bf z}^*_i:=(z_1,\ldots,z_{i-1},r_{1i},z_{i+1}\ldots,z_M)$ and $r_{1i}$ is the root of
$P_i(r)=-(\lambda_1(1-z_1)+1/b_i+\alpha_i)r+z_{i+1}/b_i$ such that $|r_{1i}|<1$.
\\
\newline
{\bf Time-limited}: according to this discipline the server continues the service of a
queue until certain exponentially distributed time of rate $\alpha_1$ will elapse or the queue becomes empty, whichever occurs first.
Consider a server visit to $Q_1$ following the same approach in Section~\ref{sec:time-limit} we find that
\begin{equation}
\mathbb{E}\big[{\bf z}^{{\bf N}^e_1}\big]=\Big(1+\frac{\alpha_1z_1}{P_1(z_1)}\Big)\mathbb{E}\big[({\bf z}_1^*)^{{\bf N}^b_1}\big]-\frac{\alpha_1z_1}{P_1(z_1)}\mathbb{E}\big[{\bf z}^{{\bf N}^b_1}\big],
\end{equation}
where ${\bf z}^*_1:=(r_{11},z_2,\ldots,z_M)$ and  $r_{11}$ is the root of $P_1(r)=\lambda_1r^2-(\lambda_1+1/b_1+\alpha_1)r+z_2/b_1$ such that $|r_{11}|<1$.  To relate $\mathbb{E}\big[{\bf z}^{{\bf N}^e_i}\big]$ to $\mathbb{E}\big[{\bf z}^{{\bf N}^b_i}\big]$ for a server
visit to $Q_i$, $i>1$, we find that
\begin{equation}
\mathbb{E}\big[{\bf z}^{{\bf N}^e_1}\big]=\Big(1+\frac{\alpha_iz_i}{P_i(z_i)}\Big)\mathbb{E}\big[({\bf z}_i^*)^{{\bf N}^b_i}\big]-\frac{\alpha_iz_i}{P_i(z_i)}\mathbb{E}\big[{\bf z}^{{\bf N}^b_i}\big],
\end{equation}
where ${\bf z}^*_i:=(z_1,\ldots,z_{i-1},r_{1i},z_{i+1}\ldots,z_M)$ and $r_{1i}$ is the root of
$P_i(r)=-(\lambda_1(1-z_1)+1/b_i+\alpha_i)r+z_{i+1}/b_i$ such that $|r_{1i}|<1$.

\section{Discussion and Conclusion}
\label{sec:conc}

In this paper, we developed a general framework to analyze polling
systems with the autonomous-server, the time-limited, and the
k-limited service discipline. The analysis of
these disciplines is based on the key idea of relating directly the
joint queue-length distribution at the beginning and the end of a server
visit.  In order to do so, we used the theory of absorbing Markov
chain. The analysis presented in this paper is restricted to the
case of service requirement with exponential distribution. We emphasize
can be extended to more general distribution such as the phase-type distributions.
For instance, Eq.~\ref{eq:ztran1} holds in the case of phase-type distribution,
however, the matrix ${\bf A_2}$ becomes a block matrix which is difficult to
invert in closed-form.

In this paper we showed that our tool is not restricted only to the disciplines that do not
verify the branching property. For example, we analyzed the exhaustive discipline. Moreover,
we claim that with an extra effort one can analyze the gated discipline for which there
already exist results in the literature.

\end{document}